\documentclass[14pt,a4paper]{article}
\RequirePackage[top=2.5cm, left=2.5cm, right=2.5cm, bottom=2.5cm]{geometry}
\usepackage[cp1251]{inputenc}
\usepackage[russian]{babel}
\usepackage{amsmath}
\usepackage{amsfonts}
\usepackage{amssymb}
\usepackage{graphicx}

\newtheorem{theorem}{Теорема}

\newtheorem{lemma}{Лемма}

\newtheorem{definition}{Определение}

\sloppypar

\begin{document}
	
УДК 517.95	

\begin{center}
\textbf{\large ANALOGUE  OF  THE  TRICOMI PROBLEM FOR THE  MIXED- TYPE  EQUATION  WITH FRACTIONAL DERIVATIVE.INVERSE PROBLEMS}

\end{center}

\begin{center}
	\textbf{R. R. Ashurov$^\textbf{1,2}$, R.T. Zunnunov$^\textbf{1}$}  
\end{center}

\begin{abstract}
In this work, an analogue of the Tricomi problem for equations of mixed type with a fractional derivative is investigated. In one part of the domain, the considered equation is a subdiffusion equation with a fractional derivative of order $\alpha\in (0,1)$ in the sense of Riemann-Liouville, and in the other it is a wave equation. Assuming the parameter $\alpha$ to be unknown, the corresponding inverse problem is studied . It was found an additional condition, that provides not only uniqueness but also existance of the desired parameter. It should be noted that the inverse problem of determining the fractional derivative for the subdiffusion  and wave equations has been studied by many mathematicians. But in the case of the Tricomi problem for a mixed-type equation, the questions of determining the fractional time derivative are studied for the first time.

\end{abstract}
\begin{center}
\textbf{\large АНАЛОГ ЗАДАЧИ ТРИКОМИ ДЛЯ УРАВНЕНИЯ СМЕШАННОГО ТИПА С ДРОБНОЙ ПРОИЗВОДНОЙ.ОБРАТНЫЕ ЗАДАЧИ}	
\end{center}	
 
\begin{center}
\textbf{Р. Р. Ашуров$^\textbf{1,2}$, Р.Т. Зуннунов$^\textbf{1}$}  
\end{center}

\begin{abstract}
	В работе исследован аналог задачи Трикоми для уравнения
	смешанного типа с дробной производной. В одной части области
	рассмотренное уравнение является уравнением субдиффузии с дробной
	производной порядка $\alpha\in (0,1)$ в смысле Римана-Лиувилля, а
	в другой - волновое уравнение. Считая параметр $\alpha$
	неизвестным, изучена соответствующая обратная задача и найдено дополнительное условие, которое обеспечивает однозначное определение  искомого параметра. 
	
\end{abstract}

\textbf{Введение.}  Известно, что для моделирования динамики сложных случайных
процессов, возникающих при аномальной диффузии в различных
областях современной науки и техники, необходимо использовать
дифференциальные уравнения дробного порядка (см., например, в
физике [1-3], в экономике [4,5], в гидрологии [6], в биологии [7]). Это
инициировало активное исследование специалистов в области
дифференциальных и псевдодифференциальных уравнений дробного
порядка (см., например,[8-14]). Имеются также
ряд работ, посвященных теории уравнений смешанного типа с участием
производных дробного порядка (см. например, [15 -22]).

Обычно обратными задачами в теории уравнений с частными
производными (целого и дробного порядка) называют такие задачи, в
которых вместе с решением дифференциального уравнения требуется
определить также тот или иной коэффициент(ы) самого уравнения,
либо правую часть, либо и коэффициент(ы), и правую часть. Интерес
к исследованию таких обратных задач обусловлен важностью их
приложений в различных разделах механики, сейсмологии, медицинской
томографии и геофизики (см. например, монографии М. М. Лаврентьев
и др.[23], С. И. Кабанихин [24,25], для
уравнений смешанного типа см. [26, 27].

В последние несколько лет интенсивно исследуется сравнительно
новый тип обратных задач, который возникает только при
рассмотрении уравнений дробного порядка. А именно, обратные задачи
по определению порядка дробной производной в различных
дифференциальных уравнениях  (см. обзорную работу [28] и
литературу там и недавние работы [34-37]). Интерес к
такой обратной задаче обясняется тем, что при рассмотрении
уравнения дробного порядка в качестве модельного уравнения для
анализа различных процессов, порядок дробной производной часто
неизвестен и его трудно измерить напрямую. Нахождение этого
неизвестного параметра интересен не только теоретически, но и
необходим для решения начально-краевых задач и изучения свойств
решений.

Данная работа посвящена изучению именно таких обратных задач для
уравнений смешанного типа.

\textbf{Постановки задач.}  Переходим к точным описаниям объектов исследования.

Дробная производная в смысле Римана-Лиувилля порядка $\alpha$,
$k-1<\alpha\leq k$, $k\in \mathbb{N}$,  от функции $h$,
определенной на $[0, \infty)$, определяется по формуле [12,c.14]
$$
\partial_t^\alpha h(t)= \frac{1}{\Gamma(k-\alpha)}\frac{d^k}{dt^k}
\int_0^t \frac{h(\tau)d\tau}{(t-\tau)^{\alpha+1-k}}, \quad t>0,
$$
при условии, что правая часть равенства существует. Здесь
$\Gamma(t)$ - гамма функция Эйлера.

Если $\alpha=k$, то дробная производная совпадает с обычным
классическим производным: $\partial_{t}^k h(t) =\frac{d^k}{dt^k}
h(t)$.

Пусть $0<\alpha < 1$. Дробный интеграл в смысле Римана-Лиувилля
порядка $\alpha -1<0$ от функции $h$ имеет вид [12, c.14]
$$
\partial_t^{\alpha-1} h(t)=\frac{1}{\Gamma
	(1-\alpha)}\int\limits_0^t\frac{h(\xi)}{(t-\xi)^{\alpha}} d\xi, \quad t>0.
$$
Отметим следующее свойство этих интегралов [12, c.104]
\begin{equation} \label{property}
\lim\limits_{t\rightarrow +0}\partial_t^{\alpha-1} h(t) = \Gamma
(\alpha)\lim\limits_{t\rightarrow +0} t^{1-\alpha} h(t).
\end{equation}

Рассмотрим уравнение, имеющее смешанный тип,
\begin{equation} \label{eq1}
f(x, t)=   \begin{cases}
\partial_t^{\alpha}u-u_{xx}, \quad t>0, \\
u_{tt} - u_{xx}, \quad t<0,
\end{cases}
\end{equation}
в двух разных областях: $G_\infty$ и $\Omega$. Первая из них -
неограниченная и имеет вид $G_\infty=G^+ \cup G^-$, при этом
$G^+=\{(x,t): -\infty<x<+\infty,\, t>0\}$ а $G^-$ - область,
находящийся в нижней полуплоскости ($t<0$) и ограниченная
характеристиками $AC: x+t=0$ и $x-t=1$ и отрезком $[0,1]$ прямой
$t=0$. А вторую (ограниченную) область определим следующим
образом: $\Omega=\Omega^+\cup \Omega^-$, где $\Omega^+=\{(x,t):
0<x<1, \, t>0\}$ и $\Omega^-=G^-$.

Сформулируем аналог задачи Трикоми для уравнения (\ref{eq1}) в
области $G_\infty$  [16,22].

\textbf{Задача} $T_\infty$. \emph{Найти решение $u(x,t)$ уравнения
	(\ref{eq1}) в области $G_\infty$, удовлетворяющее краевым
	условиям}
\[
\lim\limits_{t\rightarrow +0} t^{1-\alpha} u(x,t)=0, \,\,
-\infty<x\leq 0, \,\, 1\leq x<\infty,
\]
\begin{equation}\label{bcon}
u(x/2,-x/2)=\psi(x), \,\, 0\leq  x\leq 1.
\end{equation}
\emph{На линии $t=0$ выполняются условия склеивания}
\begin{equation}\label{gl1}
\lim\limits_{t\rightarrow +0} t^{1-\alpha}
u(x,t)=\lim\limits_{t\rightarrow -0} u(x,t),
\end{equation}
\begin{equation}\label{gl2}
\lim\limits_{t\rightarrow +0} t^{1-\alpha}
(t^{1-\alpha}u(x,t))_t=\lim\limits_{t\rightarrow -0} u_t(x,t),
\end{equation}
\emph{где} $0<x<1$.

Если $\alpha =1$, то данная задача совпадает с формулировкой
классической задачи Трикоми для уравнения смешанного
парабола-гиперболического типа. В частности, последние условия
означают непрерывность решения $u(x, t)$ и его производной по $t$
на линии изменения типа уравнения $t=0$.

Аналог задачи Трикоми для уравнения (\ref{eq1}) в области $\Omega$
имеет вид:

\textbf{Задача} $T$. \emph{Найти решение $u(x,t)$ уравнения
	(\ref{eq1}) в области $\Omega$, удовлетворяющее краевым условиям}
\begin{equation}\label{1con}
u(0,t)=u(1,t) =0, \,\, t\geq 0,
\end{equation}
и (\ref{bcon}). \emph{На линии $t=0$ выполняются те же  условия
	склеивания (\ref{gl1}) и (\ref{gl2})}.

\begin{definition}\label{defT}. Решением задачи   $T$ будем
	называть функцию $u(x,t)$, удовлетворяющую условиям задачи и
	имеющую следующую гладкость	\begin{enumerate}
		\item
		$t^{1-\alpha} u(x,t)\in C(\overline{\Omega^+})$,
		\item$t^{1-\alpha} (t^{1-\alpha}u(x,t))_t\in C(\Omega^+\cup \{(x,t): 0<x<1, \, t=0\})$,
		\item $\partial_t^\alpha u(x,t)\in C(\Omega^+)$,
		\item
		$u_{xx}(x,t)\in C(\Omega^+\cup \Omega^-),$
		\item $u_{tt}(x,t)\in C(\Omega^-)$,
		\item $u(x, t)\in \overline{\Omega^-}.$
	\end{enumerate}
\end{definition}

Задачи $T_\infty$ и $T$ также будем называть \emph{прямыми задачами.}

Основная цель настоящей работы заключается в решение \emph{обратных задач} по определению порядка дробной производной
$\alpha$ в уравнении (\ref{eq1}). Для этого сначала докажем существование и единственность и найдем решения прямых задач (в параграфе 2).  При решение обратных задач (в параграфе 3) будем
предполагать, что $f(x,t)=f(x)$ (т.е. зависит только от $x$), в
противном случае доказательства становятся технически значительно сложными.

Отметим, что прямые задачи $T_\infty$ и $ T $ для однородного
уравнения (\ref{eq1}) (т.е. при $f(x, t)\equiv 0$) впервые былы
поставлены и изучены в работах С.Х. Геккиевой [16] и К. У.
Хубиева [22] соответственно. Однако уравнение,
рассматриваемое в настоящей работе неоднородное и вид решений,
полученных в работах [16,22], не является удобным
для решения обратных задач. Поэтому мы получим приемлимый для
наших дальнейших исследований  вид решений прямых задач $T_\infty$
и $ T $. При этом подробно рассмотрим задачу $T$, а задача
$T_\infty$ изучается аналогично.

Для нахождения неизвестного порядка производной естественно надо
наложить дополнительное условие на решение $u(x,t)$. В параграфе 3
найдены такие допольнительные условия, которые гарантируют 
существование и единственность параметра $\alpha$.

\textbf{Существование и единственность решений задач Трикоми.}   Рассмотрим сначала задачу $T$.
\begin{theorem}\label{Tf}
Пусть $\psi(x)\in C[0,1]\cap C^2(0,1)$,  а функция $f(x,t)\in
C(\overline{\Omega})$ и при каждом $t$ как функция от $x$
принадлежит классу Гельдера $ C^{a}(0,1)$, $a>1/2$, и пусть $f(0,t)=f(1, t)=0$. Тогда решение задачи $T$ существует и оно
	единственно.
	
\end{theorem}

\textbf{Доказательство.}  Обозначим
\begin{equation}\label{tau}
\tau(x)=\lim\limits_{t\rightarrow +0} t^{1-\alpha} u(x,t), \,\,	0\leq x\leq 1.
\end{equation}
В силу граничных условий (\ref{1con}), имеем
	\begin{equation}\label{tau0}
	\tau(0)=\tau(1)=0.
	\end{equation}
	Введем следующий формальный ряд
	\begin{equation}\label{omega+}
	u(x,t)= \sum\limits_{k=1}^\infty \big[\Gamma(\alpha)t^{\alpha-1}
	E_{\alpha, \alpha}(-k^2 t^{\alpha}) \tau_k + \int\limits_0^t
	\eta^{\alpha-1} E_{\alpha, \alpha} (-k^2 \eta^\alpha)f_k(\eta)
	d\eta \big]\sin(k\pi x),
	\end{equation}
	где $\tau_k$ и $f_k(t)$ - коэффициенты Фурье функций $\tau(x)$ и
	$f(x, t)$ по системе $\{\sin(k\pi x)\}$ соответственно, $E_{\rho,
		\mu}(t)$ - функция Миттаг-Леффлера:
	$$
	E_{\rho, \mu}(t)= \sum\limits_{n=0}^\infty \frac{t^n}{\Gamma(\rho
		n+\mu)}.
	$$
	
	Пусть $\tau(x)$ со свойствами (\ref{tau0}) является заданной
	функцией. Тогда справедливо следующее утверждение.
	
	\begin{lemma}\label{T+f}.
		Пусть $f(x,t)$ удовлетворяет условиям Теоремы \ref{Tf}, $\tau\in
		C[0,1]\cap C^a(0,1)$, $a>1/2$ и имеют место равенства
		(\ref{tau0}). Тогда функция (\ref{omega+}) является единственным
		решением уравнения (\ref{eq1}), удовлетворяющим условиям
		(\ref{1con}), (\ref{tau}) и такое, что $t^{1-\alpha} u(x,t)\in
		C(\overline{\Omega^+})$, $\partial_t^\alpha u(x,t), \,
		u_{xx}(x,t)\in C(G)$.
	\end{lemma}

		Единственность решения сформулированной в лемме задачи доказана в
		работе  \cite{Xub}.
		
		Учитывая равенство (\ref{property}),  нетрудно проверить, что
		функция (\ref{omega+}) формально является решением уравнения
		(\ref{eq1}) и удовлетворяет условиям (\ref{1con}), (\ref{tau})
		[12, c.16] и [38]. Остается
		обосновать почленное дифференцирование ряда (\ref{omega+}).
		
		Пусть
		\[
		u_j(x,t)= \sum\limits_{k=1}^j \big[\Gamma(\alpha)t^{\alpha-1}
		E_{\alpha, \alpha}(-k^2 t^{\alpha}) \tau_k + \int\limits_0^t
		\eta^{\alpha-1} E_{\alpha, \alpha} (-k^2 \eta^\alpha)f_k(\eta)
		d\eta  \big]\sin(k\pi x).
		\]
		Тогда
		\begin{equation}\label{omega+j1}
		\frac{\partial^2}{\partial x^2}u_j(x,t)= -\sum\limits_{k=1}^j
		\big[\Gamma(\alpha)t^{\alpha-1} E_{\alpha, \alpha}(-k^2
		t^{\alpha}) \tau_k + \int\limits_0^t \eta^{\alpha-1} E_{\alpha,
			\alpha} (-k^2 \eta^\alpha)f_k(\eta) d\eta \big](k \pi)^2\sin(k\pi
		x).
		\end{equation}
		
		Функция Миттаг-Леффлера отрицательного аргумента имеет оценку
		$|E_{\rho, \mu}(-t)| \leq C(1+t)^{-1}$  [12]),
		т.е.
		\begin{equation}\label{ML}
		|E_{\alpha, \mu}(-t)|\leq C, \,\, \text{если}\,\,t<1, \quad
		|E_{\alpha, \mu}(-t)|\leq \frac{C}{t},\,\, \text{если}\,\, t\geq
		1.
		\end{equation}
		Применяя эти оценки, получим
		\begin{equation}\label{deru}
		\bigg|\frac{\partial^2}{\partial x^2}u_j(x,t)\bigg|\leq
		C\sum\limits_{k=1}^j \bigg[t^{-1}|\tau_k| +\big(|\ln
		t|+\frac{1}{\alpha} (\ln k +1)\big)|f_k(t)| \bigg].
		\end{equation}
		
		В теореме (3.1) работы А. Зигмунда [39,c.384]  доказано,
		что для произвольной функции $g(x)\in C[0,1]$, со свойствами
		$g(0)=g(1)=0$, имеет место оценка
		\[
		\sum\limits_{k=2^{n-1}+1}^{2^n} |g_k|^2 \leq
		\omega^2\bigg(\frac{1}{2^{n+1}}\bigg),
		\]
		где $\omega(\delta)$ - модуль непрерывности функции $g(x)$.
		Следовательно, по неравенству Коши - Буняковского
		\[
		\sum\limits_{k=2^{n-1}+1}^{2^n} (\ln k +1)|g_k| \leq
		\bigg(\sum\limits_{k=2^{n-1}+1}^{2^n}
		|g_k|^2\bigg)^{\frac{1}{2}}\bigg(\sum\limits_{k=2^{n-1}+1}^{2^n}
		(\ln k+1)^2\bigg)^{\frac{1}{2}}\leq C
		2^{\frac{n}{2}}n^{\frac{1}{2}}\omega\bigg(\frac{1}{2^{n+1}}\bigg),
		\]
		и наконец
		\[
		\sum\limits_{k=2}^{\infty} (\ln k +1)|g_k|=
		\sum\limits_{n=1}^\infty\sum\limits_{k=2^{n-1}+1}^{2^n} (\ln k
		+1)|g_k| \leq  C \sum\limits_{n=1}^\infty
		2^{\frac{n}{2}}n^{\frac{1}{2}}\omega\bigg(\frac{1}{2^{n+1}}\bigg).
		\]
		Очевидно, если $\omega(\delta)\leq C \delta^a$, $a>1/2$, то
		последний ряд сходится.

		Таким образом, если $\tau(x)$ и $f(x,t)$  удовлетворяют условиям
		леммы, то сумма в правой части неравенства (\ref{deru}) сходится
		при $j\rightarrow \infty$ для всех $t>0$. Отсюда следует, что
		$u_{xx}(x,t)\in C(\Omega^+)$.
		
		Согласно уравнению (\ref{eq1})
		\[
		\partial_t^\alpha u_j(x,t)=\frac{\partial^2}{\partial
			x^2}u_j(x,t)+\sum\limits_{k=1}^j f_k (t)\sin(k\pi x).
		\]
		Следовательно, из доказанного выше вытекает, что
		$\partial^\alpha_t u(x,t)\in C(\Omega^+)$.
		
		Непрерывность функции $t^{1-\alpha} u(x,t)$ в замкнутой области
		$\overline{\Omega^+}$ также вытекает из оценки функции
		Миттаг-Леффлера (\ref{ML}) и условий леммы.

	Введем обозначение
	\[
	\nu(x)=\lim\limits_{t\rightarrow +0}
	t^{1-\alpha}\big(t^{1-\alpha}u(x,t)\big)_t, \,\, 0\leq x\leq 1.
	\]
	Из граничных условий (\ref{1con}) следует
	\begin{equation}\label{nu01}
	\nu(0)=\nu(1)=0.
	\end{equation}
	Далее будем предполагать, что
	\begin{equation}\label{nua}
	\nu(x)\in C^a(0,1), \,\, a>\frac{1}{2}.
	\end{equation}

	В работе [16] доказано соотношение
	\[
	\lim\limits_{t\rightarrow +0}
	t^{1-\alpha}\big(t^{1-\alpha}u(x,t)\big)_t=\frac{1}{\Gamma(1+\alpha)}\lim\limits_{t\rightarrow
		+0} t^{1-\alpha}\big(\partial_t^{\alpha-1} u(x,t)\big)_t.
	\]
	Следовательно, в силу равенств $\big(\partial_t^{\alpha-1}
	u(x,t)\big)_t=\partial_t^\alpha u(x,t)$ и $\partial_t^\alpha
	u(x,t)=u_{xx}(x,t) + f(x,t)$, получаем
	\[
	\nu(x)=\frac{1}{\Gamma(1+\alpha)}\lim\limits_{t\rightarrow +0}
	t^{1-\alpha}\big(u_{xx}(x,t)+f(x,t)\big).
	\]
	Прежде чем вычислить этот предел отметим, что при выполнении условий леммы
	сумма (\ref{omega+j1}), умноженная на $t^{1-\alpha}$, сходится
	равномерно по $t\geq 0$. Следовательно, можно перейти в предел при
	$t\rightarrow +0$. Так как $\alpha<1$ и $E_{\alpha,
		\alpha}(0)=1/\Gamma (\alpha)$, то из равенства (\ref{omega+j1})
	будем иметь
	\[
	\nu(x)=-\frac{1}{\Gamma(1+\alpha)}\sum\limits_{k=1}^\infty
	(k\pi)^2 \tau_k \sin(k\pi
	x)=\frac{1}{\Gamma(1+\alpha)}\sum\limits_{k=1}^\infty (\tau'')_k
	\sin(k\pi x),
	\]
	где $(\tau'')_k$ - коэффициенты Фурье функции $\tau''(x)$. Эти
	ряды сходятся абсолютно и равномерно, поскольку они являются
	рядами Фурье функции $\nu(x)$ со свойствами (\ref{nu01}) и
	(\ref{nua})  [39,c.384]. Таким образом получили
	\begin{equation}\label{taunu}
	\nu(x)=\frac{1}{\Gamma(1+\alpha)} \,\tau''(x)
	\end{equation}
	- первое функциональное соотношение между $\tau(x)$ и $\nu(x)$.
	
	Предположим, что $f(x, t)$ удовлетворяет условиям теоремы и
	\begin{equation}\label{taunu12}
	\tau(x)\in C^2(0,1), \quad \nu(x) \in C^1(0,1).
	\end{equation}
	Рассмотрим задачу Коши для уравнения (\ref{eq1}) с данными
	$\tau(x)$ и $\nu(x)$. Хорошо известно, что единственное решение
	этой задачи в области $\Omega^-$, т.е. при $t<0$, дается формулой
	Даламбера
	\[
	u(x,t)=\frac{\tau(x-t)+\tau(x+t)}{2}
	+\frac{1}{2}\int\limits_{x-t}^{x+t} \nu(\xi) d\xi +
	\frac{1}{2}\int\limits_0^t \int\limits_{x-(t-\eta)}^{x+(t-\eta)}
	f(\xi, \eta) d\xi d\eta.
	\]
	Используя условие (\ref{bcon}), будем иметь
	\[
	u\bigg(\frac{x}{2},-\frac{x}{2}\bigg)=\frac{\tau (x)}{2}
	-\frac{1}{2}\int\limits_{0}^{x} \nu(\xi) d\xi + \frac{1}{2}F(x) =
	\psi(x), \quad 0\leq x\leq 1,
	\]
	где
	\[
	F(x)=\int\limits_{-\frac{x}{2}}^0\int\limits_{-\eta}^{x+\eta}f(\xi,
	\eta) d\xi d\eta.
	\]
	Отсюда нетрудно получить второе функциональное соотношение между
	$\tau(x)$ и $\nu(x)$:
	\[
	\tau'(x) - \nu(x) = 2\psi'(x) - F'(x).
	\]
	Исключая $\nu(x)$   из двух функциональных соотношений, получим
	\[
	\tau'(x) - \frac{1}{\gamma} \,\tau''(x) = 2\psi'(x) - F'(x), \quad
	\gamma=\gamma(\alpha)=\Gamma(1+\alpha).
	\]
	Применяя метод вариации постоянных, будем иметь
	\begin{equation}\label{tau1}
	\tau(x)= C_1+C_2 e^{\gamma x} +  g(x) - e^{\gamma
		x}\int\limits_0^x g'(\xi) e^{-\gamma \xi} d\xi, \quad
	g(x)=2\psi(x) - F(x),
	\end{equation}
	где $C_1$ и $C_2$   - произвольные постоянные.
	
	Из равенства (\ref{taunu}) следует
	\[
	\nu(x)= C_2 \gamma e^{\gamma x}  - g'(x)  - \gamma e^{\gamma
		x}\int\limits_0^x g'(\xi) e^{-\gamma \xi} d\xi.
	\]
	
	С помошью условия (\ref{tau0}) на функцию $\tau(x)$, определим
	константы $C_1$ и $C_2$ в равенстве (\ref{tau1}):
	\[
	C_1=-C_2, \quad C_2=\frac{1}{e^\gamma-1}\bigg( \gamma
	e^{\gamma}\int\limits_0^1 g'(\xi) e^{-\gamma \xi} d\xi -
	g(1)\bigg).
	\]
	Заметим, что $\gamma(\alpha)> 1/2$ для всех $\alpha \in (0,1)$ и
	поэтому $e^\gamma-1>0$.
	
	Очевидно, что если $f(x, t)$ и $\psi(x)$ удовлетворяют условиям
	теоремы, то функции  $\tau(x)$ и $\nu(x)$ удовлетворяют условиям
	(\ref{taunu12}) (следовательно, $\tau(x)$ также удовлетворяет
	условиям леммы, а  $\nu(x)$ - условию (\ref{nua})).
	
	Остается показать, что $t^{1-\alpha} (t^{1-\alpha}u(x,t))_t\in
	C(\Omega^+\cup \{(x,t): 0<x<1, \, t=0\})$ (см. определение \ref{defT}). Рассмотрим сначала случай $t>0$.
	
	Имеет место следующая формула дифференцирования функций
	Миттг-Леффлера (см.[10], формулу (4.3.1))
	\[
	\frac{d}{dt}\bigg[t^{\mu-1} E_{\alpha, \mu}
	(t^\alpha)\bigg]=t^{\mu-2}E_{\alpha, \mu-1} (t^\alpha).
	\]
	Отсюда получим
	\[
	\frac{d}{dt}\bigg[t^{\alpha-1} E_{\alpha, \alpha}
	(t^\alpha)\bigg]=t^{\alpha-2}E_{\alpha, \alpha-1} (t^\alpha),
	\quad \frac{d}{dt}\bigg[t^{\alpha} E_{\alpha, \alpha+1}
	(t^\alpha)\bigg]=t^{\alpha-1}E_{\alpha, \alpha} (t^\alpha).
	\]
	Следовательно,
	\[
	\frac{d}{dt}\bigg[t^{1-\alpha}\big(t^{\alpha-1} E_{\alpha, \alpha}
	(-k^2 t^\alpha)\big)\bigg]=\frac{1}{t}\bigg[(1-\alpha)E_{\alpha,
		\alpha} (-k^2 t^\alpha)+E_{\alpha, \alpha-1} (-k^2 t^\alpha)\bigg]
	\]
	и поэтому
	\[
	\big|t^{1-\alpha} (t^{1-\alpha}u_j(x,t))_t\big|\leq
	C\sum\limits_{k=1}^j \bigg[t^{-1}|\tau_k| +|f_k(t)| \bigg].
	\]
	Как было отмечено выше, данный ряд сходится при любых $t>0$.
	
	С другой стороны, $\lim\limits_{t\rightarrow +0}
	t^{1-\alpha}\big(t^{1-\alpha}u(x,t)\big)_t=\nu(x)$ и   $\nu(x)$
	удовлетворяет условию (\ref{nua})).

	Таким образом, теорема 1 полностью доказана.

Задача $T_\infty$ рассматривается совершенно аналогично, при этом
ряд Фурье следует заменить на интеграл Фурье и коэффициенты Фурье
- на преобразование Фурье.

Преобразование Фурье по переменной $x$ произвольной, интегрируемой
по $x\in (-\infty, \infty)$ функции $f(x, t)$, имеет вид
\[
\hat{f}(\xi, t)=\frac{1}{\sqrt{2\pi}}\int\limits_{-\infty}^\infty
f(x, t) e^{-i x\xi} d \xi.
\]

Условия на функцию $\psi (x)$ для существования и единственности
решения задачи $T_\infty$  остается такими же как и в теореме
\ref{Tf}, а для функции $f(x, t)$ следует требовать сходимость
интеграла
\[
\int\limits_{-\infty}^\infty \ln (|\xi|+1) |\hat{f}(\xi, t)|d \xi
<\infty,
\]
при каждом $t>0$. По неравенству Коши-Буняковского
\[
\int\limits_{-\infty}^\infty \ln (|\xi|+1) |\hat{f}(\xi, t)|d \xi
\leq ||f||_{L_2^a(-\infty, \infty)}
\bigg(\int\limits_{-\infty}^\infty\ln^2(|\xi|+1)(|\xi|^2+1)^{-a}
d\xi\bigg)^{\frac{1}{2}},
\]
где
\[
||f||^2_{L_2^a(-\infty, \infty)}=\int\limits_{-\infty}^\infty
(1+|\xi|^2)^a |\hat{f}(\xi, t)|^2 d\xi
\]
- норма функции $f(x, t)$ в пространстве Лиувилля $L_2^a(-\infty,
\infty)$, и если $a>1/2$, то правая часть последнего неравенства
ограничена.

Таким образом, справедлива теорема
\begin{theorem}\label{Tinftyf}
	Пусть $\psi(x)\in C[0,1]\cap C^2(0,1)$,  а функция $f(x,t)\in
	C(G_\infty)$ и при каждом $t$ как функция от $x$ принадлежит
	классу Лиувилля $L_2^a(-\infty, \infty)$, $a>1/2$. Тогда решение
	задачи $T_\infty$ единственно и оно имеет вид
	\begin{equation}\label{omega-}
	u(x,t)= \frac{1}{\sqrt{2\pi}}\int\limits_{-\infty}^\infty
	\big[\Gamma(\alpha)t^{\alpha-1} E_{\alpha, \alpha}(-\xi^2
	t^{\alpha}) \hat{\tau}(\xi) + \int\limits_0^t \eta^{\alpha-1}
	E_{\alpha, \alpha} (-\xi^2 \eta^\alpha)\hat{f}(\xi, \eta) d\eta
	\big]e^{ix\xi} d \xi.
	\end{equation}
	
\end{theorem}

\textbf{Обратные задачи по определению порядка производных.}
Теперь предположим, что порядок дробной производной $\alpha$ в
уравнении (\ref{eq1}) является неизвестным и рассмотрим обратную
задачу по определению этого параметра. Исследуем сначала в случае
задачи $T$.

Как было отмечено выше, в этом параграфе будем считать, что
$f(x,t)=f(x)$. При этом естественно предполагать, что $f(x)$ и
$\tau(x)$ не являются одновременно тождественными нулями, так как
в противном случае, в силу леммы \ref{T+f}, решение прямой задачи
$u(x,t)\equiv 0$ в области $\Omega^+$ при любых $\alpha$, что не
позволяет определить этот параметр.

Определение правильного порядка дробной производной  играет важную
роль в моделирование различных процессов. Соответствующая обратная
задача для уравнений субдиффузии рассматривалась рядом авторов
(см. обзорную статью Li, Liu, Yamamoto [28] и ссылки в
нем, [29 - 37]). Отметим, что во всех известных
работах уравнение субдиффузии рассматривалось в ограниченной
области $\Omega \subset \mathbb{R}^N $. Следует также отметить,
что в работах [29 - 32] следующее соотношение:
\begin{equation}\label{ex1}
u(х_0, t) = h (t), \, \, 0 <t <T,
\end{equation}
в точке мониторинга $ x_0\in\overline{\Omega} $ принималось как
дополнительное условие. Но это условие, как правило (исключение
составляет работа J. Janno [32], где доказаны и
единственность, и существование), может обеспечить только
единственность решения обратной задачи [29 -31].
Авторы статьи Ашуров и Умаров \cite{AU} рассмотрели значение
проекции решения на первую собственную функцию эллиптической части
уравнение субдиффузии в качестве дополнительной информации.
Следует отметить, что результаты работы [34] применимы,
только если первое собственное значение равно нулю. Без этого
ограничения существование и единственность неизвестного порядка
дробной производной были доказаны в недавней работе Алимова и
Ашурова [35]. В этой статье дополнительное условие имеет вид
$ || u (x, t_0) ||^2 = d_0 $ и граничное условие не обязательно
однородное. В работах [33] и [36] авторы исследовали
обратную задачу для одновременного определения порядка дробной
производной по времени и функции источника в уравнениях
субдиффузии.

В настоящей работе в случае задачи $T$ для определения
неизвестного параметра $\alpha$ зададим допольнительное условие на решения в виде:
\begin{equation}\label{excon1}
\int\limits_0^1 u(x,t_0) \sin(k_0\pi x) dx=d_0,
\end{equation}
где $d_0$ - произвольное заданное число, $t_0$ - положительное
число, которое определено в лемме \ref{monoton}, приведенной ниже
и $k_0\geq 1$ - произвольное целое такое, что $f^2_{k_0} +
\tau^2_{k_0} \neq 0$ (очевидно, такие числа существуют, так как
$f(x)$ и $\tau(x)$ не являются одновременно тождественными
нулями).

Начально-краевую задачу $T$ вместе с дополнительным условием
(\ref{excon1}) назовем \emph{обратной задачей} по определению
параметра $\alpha$. Если $u(x,t)$ решение задачи $T$ и параметр
$\alpha$ удовлетворяет условию (\ref{excon1}), то пару
$\{u(x,t),\alpha\}$ назовем \emph{решением обратной задачи}.

Так как  $f(x,t)= f(x)$, то в силу формулы (см. [10],
формула (4.4.4))
\[
\int\limits_0^t \eta^{\alpha-1} E_{\alpha, \alpha} (-k^2
\eta^\alpha) d\eta= t^\alpha E_{\alpha, \alpha+1} (-k^2 t^\alpha)
\]
будем иметь (см. (\ref{omega+}))
\[
E(\alpha)\equiv\int\limits_0^1 u(x,t_0) \sin(k_0\pi x)
dx=2\big(\Gamma(\alpha)t_0^{\alpha-1} E_{\alpha, \alpha}(-k_0^2
t_0^\alpha)\tau_{k_0}+t_0^{\alpha} E_{\alpha, \alpha+1}(-k_0^2
t_0^\alpha)f_{k_0}\big).
\]

Всюду далее будем предполагать, что $\alpha\in [\alpha_0, 1]$,
$0<\alpha_0<1$. Тогда функция $E(\alpha)$ является непрерывно
дифференцируемой на сегменте  $\alpha\in [\alpha_0,1]$.

Перепишем допольнительное условие (\ref{excon1})  в виде
следующего функционального уравнения относительно $\alpha$:
\begin{equation}\label{excon2}
E(\alpha)=d_0.
\end{equation}
Для того, чтобы это уравнение имело решение, очевидно, что число
$d_0$ не может задаваться произвольным образом. Необходимое
условие разрешимости уравнения (\ref{excon2}) является выполнение
следующего включения
\begin{equation}\label{excon3}
d_0\in [\min\limits_{[\alpha_0,1]} E(\alpha),\,
\max\limits_{[\alpha_0,1]} E(\alpha)],
\end{equation}
т.е. число $d_0$ должно находиться в области значений функции
$E(\alpha)$. Другими словами, если условие  (\ref{excon3}) не
выполняется, то ни при каких $\alpha\in [\alpha_0, 1]$ не
существует решение $u(x, t)$ прямой задачи, которое удовлетворяет
условию (\ref{excon1}). Всюду далее в этом параграфе будем
предполагать, что это условие выполнено.

\begin{lemma}\label{monoton}. Пусть выполнены условия теоремы
	\ref{Tf}. Найдется число $T_{0}= T_{0}(k_0, \alpha_0)$ такое, что
	при $t_0 \geq T_0 $ функция $E(\alpha)$ является строго монотонной
	по $\alpha\in [\alpha_0, 1]$.
	
\end{lemma}

Из этой леммы очевидным образом вытекает следующий основной
результат относительно задачи $T$.

\begin{theorem}\label{it}. Пусть выполнены условия теоремы
	\ref{Tf} и  $t_0 \geq T_0 $. Тогда решение обратной задачи
	существует и оно единственно.
\end{theorem}

Прежде чем доказать лемму  \ref{monoton}, изучим поведение
производных по $\alpha$ функций $e_{\lambda, 1}
(\alpha)\equiv\Gamma(\alpha) t_0^{\alpha-1}E_{\alpha,
	\alpha}(-\lambda t_0^\alpha)$ и $e_{\lambda, 2} (\alpha)\equiv
t_0^\alpha E_{\alpha,\alpha+1}(-\lambda t_0^\alpha)$. Данные
функции участвуют в определении функции $E(\alpha)$.

\begin{lemma}\label{el1}.
	Существует число $T_1= T_1 (\alpha_0, \lambda)>0$ такое, что для
	всех $t_0\geq T_1$ и $\alpha\in [\alpha_0, 1]$ производная по
	$\alpha$ функции $e_{\lambda, 1} (\alpha)$ отрицательна и
	справедлива оценка
	\begin{equation}\label{el1e}
	\bigg|\frac{d}{d\alpha} e_{\lambda, 1}(\alpha)\bigg|\leq
	\frac{C\cdot \ln t_0}{\alpha_0 \lambda^2 t_0^{\alpha+1}}.
	\end{equation}
\end{lemma}

\textbf{Доказательство.} Обозначим через $ \delta (\theta) $ контур, пробегаемый по
	направлении неубывания велечины $ \arg \zeta $ и состоящий из
	следующих частей: луч $ \arg \zeta = - \theta $ с $ | \zeta | \geq
	1 $, arc $ - \theta \leq \arg \zeta \leq \theta $, $ | \zeta | = 1
	$ и луч $ \arg \zeta = \theta $, $ | \zeta | \geq 1 $. Если $ 0
	<\theta <\pi $, то контур $  \delta (\theta) $ делит всю
	комплексную $ \zeta $ -плоскость на две неограниченные части, а
	именно $ G^{(-)} (\theta) $ слева от $ \delta (\theta) $ и $
	G^{(+)} (\theta) $ справа от него. Контур $ \delta (\theta) $
	называется путем Ханкеля [11, c.126].

	Пусть $\theta = \frac{3\pi}{4}\alpha$, $\alpha\in [\alpha_0, 1)$.
	Тогда $-\lambda t_0^\alpha\in G^{(-)}(\theta)$ и в силу
	определения контура $\delta(\theta)$, будем иметь [11, c.135]
	\begin{equation}\label{Erho}
	e_{\lambda, 1}(\alpha)= - \frac{\Gamma (\alpha)}{\lambda^2
		t_0^{\alpha+1} \Gamma (-\alpha)}+\frac{\Gamma (\alpha)}{2\pi i
		\alpha \lambda^2
		t_0^{\alpha+1}}\int\limits_{\delta(\theta)}\frac{e^{\zeta^{1/\alpha}}\zeta^{\frac{1}{\alpha}+1}}{\zeta+\lambda
		t_0^\alpha} d\zeta = p_1(\alpha)+p_2(\alpha).
	\end{equation}

	Нетрудно оценить производную $p'_1(\alpha)$. Действительно, пусть
	$\Psi(\alpha)$ логарифмическая производная гамма функции
	$\Gamma(\alpha)$ (определение и свойства функции $\Psi$ см. в
	работе [40]). Тогда $\Gamma'(\alpha) = \Gamma (\alpha)
	\Psi(\alpha)$ и поэтому
	$$
	p_1'(\alpha)=\Gamma (\alpha) \frac{\ln t_0 - \Psi
		(-\alpha)}{\lambda^2 t_0^{\alpha +1}\Gamma (-\alpha)} +
	\frac{\Gamma (\alpha) \Psi(\alpha) }{\lambda^2 t_0^{\alpha+1}
		\Gamma (-\alpha)}= p'_{11}(\alpha)+p'_{12}(\alpha).
	$$
	Так как
	$$
	\frac{1}{\Gamma(-\alpha)}=-\frac{\alpha}{\Gamma(1-\alpha)}=-\frac{\alpha(1-\alpha)}{\Gamma(2-\alpha)},
	\quad \Psi(-\alpha)=\Psi(1-\alpha)
	+\frac{1}{\alpha}=\Psi(2-\alpha)+\frac{1}{\alpha}-\frac{1}{1-\alpha},
	$$
	то функцию $p_{11}'(\rho)$  можно переписать ввиде
	\begin{equation}\label{f11}
	p_{11}'(\alpha)=\frac{\Gamma (\alpha)}{\lambda^2
		t_0^{\alpha+1}\Gamma(2-\alpha)}\cdot \big(\alpha(1-\alpha)[\Psi
	(2-\alpha)-\ln t_0]+1-2\alpha\big)=-\frac{\Gamma
		(\alpha)}{\lambda^2 t_0^{\alpha+1}\Gamma(2-\alpha)}\cdot
	p_{0}(\alpha).
	\end{equation}
	Пусть $\gamma\approx 0,57722$ - константа Эйлера-Машероны. Тогда
	$\Psi(2-\alpha)< 1-\gamma$ и поэтому
	\[
	p_{0}(\alpha)>\alpha (1-\alpha)[\ln t_0 -(1-\gamma)])+2\alpha-1.
	\]
	Если $t_0=e^{1-\gamma} e^{2/\alpha}$, то $\alpha (1-\alpha)[\ln
	t_0 -(1-\gamma)])+2\alpha-1=1$. Следовательно, выбрав $t_0\geq
	T_0\equiv e^{1-\gamma} e^{2/\alpha_0}$, будем иметь
	$p_{0}(\alpha)> 1$.

	Таким образом, согласно равенству (\ref{f11}), для всех таких
	$t_0$ следует оценка
	\[
	p_{11}'(\alpha)\leq-\frac{\Gamma(\alpha)}{\lambda^2
		t_0^{\alpha+1}}.
	\]
	
	Перепишем $p_{12}(\alpha)$ ввиде
	\begin{equation}\label{f12}
	p'_{12}(\alpha)=\frac{\Gamma (\alpha) \Psi(\alpha) }{\lambda^2
		t_0^{\alpha+1} \Gamma (-\alpha)}=\frac{\Gamma (\alpha)
		\Psi(\alpha) }{\lambda^2 t_0^{\alpha+1}}\,\,\frac{\alpha
		(1-\alpha)}{ \Gamma (2-\alpha)}=\frac{\Gamma (\alpha) }{\lambda^2
		t_0^{\alpha+1}}\,\,\frac{(\Psi(\alpha+1)\alpha-1) (1-\alpha)}{
		\Gamma (2-\alpha)}.
	\end{equation}
	Так как $-\gamma< \Psi(\alpha+1)< 1-\gamma$, то
	$p'_{12}(\alpha)<0$ для всех $\alpha\in [\alpha_0,1)$.
	Следовательно, всилу оценки $p_{11}'(\alpha)$, будем иметь
	\begin{equation}\label{f1}
	p_{1}'(\alpha)\leq-\frac{\Gamma(\alpha)}{\lambda^2
		t_0^{\alpha+1}}.
	\end{equation}

	Для того, чтобы оценить $p'_2(\alpha)$, введем обозначение
	\[
	F(\zeta, \alpha)=\frac{1}{2\pi i \alpha\lambda^2
		t_0^{\alpha+1}}\cdot
	\frac{e^{\zeta^{1/\alpha}}\zeta^{1/\alpha+1}}{\zeta+\lambda
		t_0^\alpha}.
	\]
	Заметим, что область интегрирования $\delta(\theta)$ также зависит
	от $\alpha$. Что бы это учесть при дифференцировании функции
	$p'_2(\alpha)$, перепишем интеграл (\ref{Erho}) в виде:
	\[
	\frac{1}{\Gamma(\alpha)}p_2(\alpha)=p_{2+}(\alpha)+p_{2-}(\alpha)+p_{21}(\alpha),
	\]
	где
	\[
	p_{2\pm}(\alpha)=e^{\pm i \theta}\int\limits_1^\infty F(s\,e^{\pm
		i \theta}, \alpha)\, ds,
	\]
	
	\[
	p_{21}(\alpha) = i \int\limits_{-\theta}^{\theta} F(e^{i y},
	\alpha)\, e^{iy} dy= i\theta \int\limits_{-1}^{1} F(e^{i \theta
		s}, \alpha)\, e^{i\theta s} ds.
	\]
	Рассмотрим функцию $p_{2+}(\alpha)$. Поскольку $\theta =
	\frac{3\pi}{4}\alpha$ и $\zeta= s\, e^{i\theta}$,  то
	\[
	e^{\zeta^{1/\alpha}}=e^{\frac{1}{\sqrt{2}} (i-1)\,s^{\frac{1}{\alpha}}}.
	\]
	Производная функции $ p_{2+}(\alpha)$ имеет вид
	$$
	p_{2+}'(\alpha)=I\cdot\int\limits_1^\infty \frac{e^{\frac{i-1}{\sqrt{2}}
			\,s^{1/\alpha}}s^{\frac{1}{\alpha}+1}\,e^{2ia\alpha}\big[\frac{1}{\alpha^2}(\frac{1-i}{\sqrt{2}}s^{\frac{1}{\alpha}}-1)\ln
		s+2ia -\frac{1}{\alpha}-\ln t_0-\frac{ia s\,e^{ia\alpha}+\lambda
			t_0^\alpha \ln t_0}{s\,e^{ia\alpha}+\lambda
			t_0^\alpha}\big]}{s\,e^{ia\alpha}+\lambda t_0^\alpha} ds,
	$$
	где $I=e^{ia}(2\pi i \alpha\lambda^2 t_0^{\alpha+1})^{-1}$ и
	$a=\frac{3\pi}{4}$. В силу неравенства $|s\,e^{ia\alpha}+\lambda
	t_0^\alpha|\geq \lambda t_0^\alpha$ напишем
	$$
	|p_{2+}'(\alpha)|\leq \frac{C}{\alpha \lambda^3
		t_0^{2\alpha+1}}\int\limits_1^\infty
	e^{-\frac{1}{\sqrt{2}}\,s^{1/\alpha}}s^{\frac{1}{\alpha}+1}\,\big[\frac{1}{\alpha^2}s^{1/\alpha}\ln
	s+\ln t_0\big] ds.
	$$
	
	\begin{lemma}\label{I}. \  Пусть $0<\alpha\leq 1$ и $m\in \mathbb{N}$.
		Тогда
		\[
		J(\alpha)=\frac{1}{\alpha}\int\limits_1^\infty
		e^{-\frac{1}{\sqrt{2}}s^{\frac{1}{\alpha}}} s^{\frac{m}{\alpha}+1} ds\leq
		C_m.
		\]
	\end{lemma}
	\textbf{Доказательство.} Если $r=s^{\frac{1}{\alpha}}$, то
		\[
		s=r^\alpha, \quad ds = \alpha r^{\alpha-1} dr.
		\]
		Поэтому,
		\[
		J(\alpha)=\int\limits_1^\infty e^{-\frac{1}{\sqrt{2}}r} r^{m-1+2\alpha}
		dr\leq \int\limits_1^\infty e^{-\frac{1}{\sqrt{2}}r} r^{m+1} dr = C_m.
		\]
		Лемма \ref{I} доказана.

	Применение леммы дает (заметим, что $\frac{1}{\alpha} \ln s<
	s^{\frac{1}{\alpha}}$ при $s\geq 1$)
	$$
	|p_{2+}'(\alpha)|\leq \frac{C}{\lambda^3
		t_0^{2\alpha+1}}\,\big[\frac{C_3}{\alpha}+C_1\ln t_0\big]\leq
	\frac{C}{\lambda^3 t_0^{2\alpha+1}}\,\big[\frac{1}{\alpha}+\ln
	t_0\big].
	$$
	
	Функция $p'_{2-}(\alpha)$ имеет точно такую же оценку.
	
	Теперь рассмотрим функцию $p_{21}(\alpha)$. Для производной имеем
	\[
	p'_{21}(\alpha)=\frac{a}{2\pi i \lambda^2
		t_0^{\alpha+1}}\cdot\int\limits_{-1}^1\frac{e^{e^{ias}}\,
		e^{ias}\,e^{2ia\alpha s}\big[2ias-\ln t_0-\frac{ias e^{ia\alpha
				s}+\lambda t_0^\alpha \ln t_0}{e^{ia\alpha s}+\lambda
			t_0^\alpha}\big]}{e^{ia\alpha s}+\lambda t_0^\alpha} ds.
	\]
	Поэтому
	\[
	|p'_{21}(\alpha)|\leq C\,\frac{\ln t_0}{\lambda^3
		t_0^{2\alpha+1}}.
	\]
	
	Таким образом, имеет место оценка
	\[
	\bigg|\bigg(\frac{p_2(\alpha)}{\Gamma(\alpha)}\bigg)'\bigg|=|f'_{2+}(\alpha)+f'_{2-}(\alpha)+f'_{21}(\alpha)|\leq
	\frac{C}{\lambda^3 t_0^{2\alpha+1}}\,\big[\frac{1}{\alpha}+\ln
	t_0\big].
	\]
	
	Оценивая функции $p_{2\pm}(\alpha)$ и $p_{21}(\alpha)$ как и выше,
	получим
	\[
	\bigg|\frac{p_2(\alpha)}{\Gamma(\alpha)}\bigg|\leq
	\frac{C}{\lambda^3 t_0^{2\alpha+1}}.
	\]
	Следовательно, так как
	$\Gamma'(\alpha)=\Gamma(\alpha)\Psi(\alpha)$, то
	\begin{equation}\label{f'2}
	|p'_2(\alpha)|=\bigg|\bigg(\frac{p_2(\alpha)}{\Gamma(\alpha)}\bigg)'\cdot
	\Gamma(\alpha)+\frac{p_2(\alpha)}{\Gamma(\alpha)}\cdot
	\Gamma'(\alpha)\bigg|\leq \frac{C\cdot \Gamma(\alpha)}{\lambda^3
		t_0^{2\alpha+1}}\,\big[\frac{1}{\alpha}+\ln t_0\big].
	\end{equation}
	
	Отсюда, в силу оценки (\ref{f1}) и равенства
	$\Gamma(\alpha+1)=\Gamma(\alpha) \alpha$, будем иметь
	\[
	\frac{d}{d\alpha} e_{\lambda, 1}(\alpha)\leq
	-\frac{\Gamma(\alpha)}{\lambda^2 t_0^{\alpha+1}}+C\cdot
	\Gamma(\alpha)\frac{1/\alpha+\ln t_0}{\lambda^3 t_0^{2\alpha+1}}=
	-\frac{\Gamma(\alpha+1)}{\alpha\lambda^2
		t_0^{\alpha+1}}\bigg[1-C\cdot \frac{1/\alpha+\ln t_0}{\lambda
		t_0^{\alpha}}\bigg].
	\]
	Следовательно, существует число $T_1= T_1 (\alpha_0, \lambda)>0$
	такое, что для всех $t_0\geq T_1$ и $\alpha\in [\alpha_0, 1]$
	справедлива оценка
	\[
	\frac{d}{d\alpha} e_{\lambda, 1}(\alpha)\leq
	-\frac{1}{2\alpha\lambda^2 t_0^{\alpha+1}}<0.
	\]
	С другой стороны, из равенств (\ref{f11}), (\ref{f12}) и оценки
	(\ref{f'2}), следует справедливость соотношения (\ref{el1e}) для
	всех $t_0\geq T_1$ и $\alpha\in [\alpha_0, 1]$. Лемма 3 доказана.

\begin{lemma}\label{el2}.
	Существует число $T_2= T_2 (\alpha_0, \lambda)>0$ такое, что для
	всех $t_0\geq T_2$ и $\alpha\in [\alpha_0, 1]$ производная по
	$\alpha$ функции $e_{\lambda, 2} (\alpha)$ отрицательна и
	справедлива оценка
	\begin{equation}\label{el2e}
	\frac{d}{d\alpha} e_{\lambda, 2}(\alpha)= -\frac{1}{\lambda^2
		t_0^{\alpha}}\cdot\frac{(1-\alpha)[\ln t_0 - \Psi (2-\alpha)]+1}{
		\Gamma (2-\alpha)} +\frac{ \ln t_0}{\lambda^3 t_0^{2\alpha}}\cdot
	r_2(\alpha, t_0, \lambda),
	\end{equation}
	где
	\[
	|r_2(\alpha, t_0, \lambda)|\leq C_2,
	\]
	при этом $C_2$ - абсолютная константа.
\end{lemma}

\textbf{Доказательство.} Параметр $\theta$  контура
	$\delta(\theta)$ выберем также как и выше. Тогда, по определению
	контура, имеем  [11, c.135]
	\[
	e_{\lambda, 2}(\alpha)= \frac{1}{\lambda}- \frac{1}{\lambda^2
		t_0^{\alpha} \Gamma (1-\alpha)}+\frac{1}{2\pi i \alpha \lambda^2
		t_0^{\alpha}}\int\limits_{\delta(\theta)}\frac{e^{\zeta^{1/\alpha}}\zeta}{\zeta+\lambda
		t_0^\alpha} d\zeta = \frac{1}{\lambda}- q_1(\alpha)+q_2(\alpha).
	\]

	Производная функции $q_1(\alpha)$ имеет вид
	$$
	q_1'(\alpha)=-\frac{\ln t_0 - \Psi (1-\alpha)}{\lambda^2
		t_0^\alpha \Gamma (1-\alpha)}.
	$$
	Следовательно, при $t_0\geq e$ будем иметь
	\begin{equation}\label{edq1}
	-q_1'(\alpha)=\frac{1}{\lambda^2 t_0^\alpha}\cdot
	\frac{(1-\alpha)[\ln t_0 - \Psi (2-\alpha)]+1}{ \Gamma
		(2-\alpha)}\geq \frac{1}{\lambda^2 t_0^{\alpha}}.
	\end{equation}
	
	Напомним, в силу выбора $\theta$, вдоль контура $\delta(\theta)$
	выполняются соотношения
	\[
	\zeta= s\, e^{i\theta}, \quad e^{\zeta^{1/\alpha}}=e^{\frac{1}{\sqrt{2}}
		(i-1)\,s^{\frac{1}{\alpha}}},\quad |\zeta+\lambda t_0^\alpha|\geq
	\lambda t_0^\alpha.
	\]
	Поэтому, повторяя аналогичные рассуждения, что и в доказательстве
	леммы \ref{el1}, получим оценку
	\[
	|q_2'(\alpha)|\leq C \,\frac{1/\alpha+\ln t_0}{\lambda^3 t_0
		^{2\alpha}}.
	\]
	
	Таким образом, имеем
	\[
	\frac{d}{d\alpha} e_{\lambda,
		2}(\alpha)\leq -\frac{1}{\lambda^2 t_0^{\alpha}}+C\cdot
	\frac{1/\alpha+\ln t_0}{\lambda^3 t_0^{2\alpha}}=
	-\frac{1}{\lambda^2 t_0^{\alpha}}\bigg[1-C\cdot \frac{1/\alpha+\ln
		t_0}{\lambda t_0^{\alpha}}\bigg].
	\]
	Следовательно, существует число $T_2= T_2 (\alpha_0, \lambda)>0$
	такое, что для всех $t_0\geq T_2$ и $\alpha\in [\alpha_0, 1]$
	данная производная является отрицательной. Очевидно, что при тех
	же $t_0$  имеет место оценка (\ref{el2e}) (см. соотношение
	(\ref{edq1})). Лемма 5 доказана.
	
Переходим к доказательству леммы \ref{monoton}. Для этого
достаточно показать существование такого числа $T_0$, что при
$t_0\geq T_0$ производная функции
\[
E(\alpha)= 2\, [ e_{k_0^2, 1}(\alpha) \tau_{k_0}+e_{k_0^2,
	2}(\alpha) f_{k_0}]
\]
сохраняет знак. Напомним, что число $k_0$ выбрано таким образом,
что $f_{k_0}^2+\tau_{k_0}^2\neq 0$.

Рассмотрим отдельно следующие два возможных случая: 1)
$\tau_{k_0}\cdot f_{k_0}\geq 0$ и 2) $\tau_{k_0}\cdot f_{k_0}< 0$.

Так как производные функций $e_{k_0^2, 1}(\alpha)$ и $e_{k_0^2,
	2}(\alpha)$ отрицательны при достаточно больших $t_0$, то  в
случае 1) производная $E'(\alpha)$ сохраняет знак.

Согласно лемме \ref{el2}  в случае 2) имеем
\[
E'(\alpha)= 2\,\bigg[e_{k_0^2, 1}(\alpha) \tau_{k_0} -f_{k_0}\cdot
\frac{1}{\lambda^2 t_0^{\alpha}}\cdot\frac{(1-\alpha)[\ln t_0 -
	\Psi (2-\alpha)]+1}{ \Gamma (2-\alpha)} +f_{k_0}\cdot\frac{ \ln
	t_0}{\lambda^3 t_0^{2\alpha}}\cdot r_2(\alpha, t_0,
\lambda)\bigg].
\]
Из леммы \ref{el1} и оценки функции $r_2$ следует, что при
достаточно больших  $t_0$ знак производной   $E'(\alpha)$
совпадает со знаком $-f_{k_0}$ при всех $\alpha\in [\alpha_0, 1]$.

Таким образом, лемма \ref{monoton} и следовательно теорема
\ref{it}, доказаны полностью.

\

Обратная задача в случае задачи $T_\infty$ рассматривается
совершенно аналогично. В этом случае для определения
неизвестного параметра $\alpha$ зададим допольнительное условие в виде:
\begin{equation}\label{excon4}
\hat{u}(\xi_0,t_0) =d_1,
\end{equation}
где $d_1$ - произвольное заданное число, $t_0$ - положительное
число, которое определено в лемме \ref{monoton} и $\xi_0\in
(-\infty, \infty)$  такое, что $\hat{f}^2(\xi_0) + \hat{\tau}^2
(\xi_0) \neq 0$ (очевидно, такое $\xi_0$ существуют, так как
$f(x)$ и $\tau(x)$ не являются одновременно тождественными
нулями).

Задачу $T_\infty$ с дополнительным условием (\ref{excon4}) назовем
\emph{второй обратной задачей} по определению неизвестного порядка
производной $\alpha$.

В рассматриваемом случае функция $E(\alpha)$  будет иметь вид:
\[
E(\alpha)\equiv\hat{u}(\xi_0,t_0)=\frac{1}{\sqrt{2\pi}}\big[\Gamma(\alpha)t_0^{\alpha-1}
E_{\alpha, \alpha}(-\xi_0^2
t_0^\alpha)\hat{\tau}(\xi_0)+t_0^{\alpha} E_{\alpha,
	\alpha+1}(-\xi_0^2 t_0^\alpha)\hat{f}(\xi_0)\big].
\]
Для того, чтобы вторая обратная задача имела решение, будем
считать, что число $d_1$ удовлетворяет условию (\ref{excon2}).

Сформулируем соответсвующий результат.

\begin{theorem}. Пусть выполнены условия теоремы
	\ref{Tinftyf} и  $t_0 \geq T_0 $. Тогда решение второй обратной
	задачи существует и оно единственно.
\end{theorem}

Авторы приносят глубокую благодарность Ш.А. Алимову за обсуждения
результатов работы.

\end{document}